\documentclass[11pt]{amsart}
\usepackage{amsxtra}
\usepackage{amssymb}
\usepackage{mathtools}
\theoremstyle{plain}

\newcommand{\cleqn}{\setcounter{equation}{0}}
\newcommand{\clth}{\setcounter{theorem}{0}}
\newcommand {\sectionnew}[1]{\section{#1}\cleqn\clth}

\newtheorem{theorem}{Theorem}[section]
\newtheorem{lemma}[theorem]{Lemma}
\newtheorem{definition-theorem}[theorem]{Definition-Theorem}
\newtheorem{proposition}[theorem]{Proposition}
\newtheorem{corollary}[theorem]{Corollary}
\newtheorem{definition}[theorem]{Definition}
\newtheorem{example}[theorem]{Example}
\newtheorem{remark}[theorem]{Remark}
\newtheorem{conjecture}[theorem]{Conjecture}

\newcommand \bth[1] { \begin{theorem}\label{t#1} }
\newcommand \ble[1] { \begin{lemma}\label{l#1} }

\newcommand \bpr[1] { \begin{proposition}\label{p#1} }
\newcommand \bco[1] { \begin{corollary}\label{c#1} }
\newcommand \bde[1] { \begin{definition}\label{d#1}\rm }
\newcommand \bex[1] { \begin{example}\label{e#1}\rm }
\newcommand \bre[1] { \begin{remark}\label{r#1}\rm }
\newcommand \bcj[1] { \begin{conjecture}\label{j#1}\rm }

\renewcommand {\eth} { \end{theorem} }
\newcommand {\ele} { \end{lemma} }

\newcommand {\epr} { \end{proposition} }
\newcommand {\eco} { \end{corollary} }
\newcommand {\ede} { \end{definition} }
\newcommand {\eex} { \end{example} }
\newcommand {\ere} { \end{remark} }
\newcommand {\ecj} { \end{conjecture} }

\newcommand {\enota} { \end{notation} }
\newcommand \thref[1]{Theorem \ref{t#1}}
\newcommand \leref[1]{Lemma \ref{l#1}}

\newcommand \cjref[1]{Conjecture \ref{j#1}}

\newcommand \reref[1]{Remark \ref{r#1}}
\newcommand \lb[1]{\label{#1}}

\def \Cset {{\mathbb C}}
\def \KK {{\mathbb K}}
\def \Zset {{\mathbb Z}}
\def \Nset {{\mathbb N}}

\def \B  {{\mathcal{B}}}

\def \QQ {{\mathcal{Q}}}
\def \PP {{\mathcal{P}}}

\def \UU {{\mathcal{U}}}
\def \RR {{\mathcal{R}}}

\def \EE {{\mathcal{E}}}

\def \De {\Delta}   

\def \al {\alpha}
\def \be {\beta}
\def \vpi {\varpi}
\def \la {\lambda}

\def \Sig {\Sigma}
\def \sig {\sigma}

\def \vp {\varphi}

\def \sig{\sigma}

\def \mt  {\mapsto}

\def \rcor {\rangle}
\def \lcor {\langle}
\def \del {\partial}

\def \ol {\overline}

\def \wh {\widehat}

\def \id { {\mathrm{id}} }

\def \rank { {\mathrm{rank}} }

\def \g  {\mathfrak{g}}   

\def \sl {\mathfrak{sl}} 
\def \h  {\mathfrak{h}}

\def \n  {\mathfrak{n}}

\def \sl {\mathfrak{sl}}

\DeclareMathOperator \Span { {\mathrm{Span}} }

\DeclareMathOperator \HSpec { \textit{H} \mbox{-} {\mathrm{Spec}}}

\DeclareMathOperator \Symp  { {\mathrm{Symp}} }
\DeclareMathOperator \Dix { {\mathrm{Dix}} }
\DeclareMathOperator \Prim { {\mathrm{Prim}} }

\newcommand \Spec { {\mathrm{Spec}} }
\begin{document}
\title[Separating Ore sets]
{Separating Ore sets for prime ideals \\ 
of quantum algebras}
\author[Si\^an Fryer]{Si\^an Fryer}
\address{School of Mathematics \\ University of Leeds \\ Leeds LS2 9JT\\ UK.}
\email{s.fryer@leeds.ac.uk}
\author[Milen Yakimov]{Milen Yakimov}
\address{
Department of Mathematics \\
Louisiana State University \\
Baton Rouge, LA 70803
U.S.A.
}
\email{yakimov@math.lsu.edu}
\thanks{The research of S.F. was supported by an EPSRC Doctoral Prize
Fellowship and that of M.Y. by the NSF grant DMS-1303038, a Scheme 2
grant from the London Mathematical Society and Louisiana
Board of Regents grant Pfund-403.}
\begin{abstract} Brown and Goodearl stated a conjecture that provides an explicit description of the topology of the spectra of quantum algebras.
The conjecture takes on a more explicit form if there exist separating Ore sets for all incident pairs of torus invariant prime ideals of the given 
algebra. We prove that this is the case  
for the two largest classes of algebras of finite Gelfand--Kirillov dimension that fit the setting of the conjecture: the quantized coordinate rings of all simple 
algebraic groups and the quantum Schubert cell algebras for all symmetrizable Kac--Moody algebras. 
\end{abstract}
\maketitle
\sectionnew{The Brown--Goodearl conjecture and overview of results in the paper}
\lb{BGcnj}
\subsection{Topology of spectra on noncommutative algebras}
\label{1.1}
The study of the topology of the spectra of noncommutative algebras has been a major direction in ring theory, starting 
with the Dixmier program \cite{Dix} on the spectra of universal enveloping algebras. It is related to problems in Lie theory, 
Poisson geometry and quantization. The description of the topology of spectra is equivalent to classifying the inclusions between 
prime ideals of the given algebra.

In the early 90's Hodges--Levasseur \cite{HL,HLT} and Joseph \cite{J0,J} initiated the study of the spectra of quantum function algebras 
on complex simple groups $R_q[G]$. Since then many results have been obtained on the structure of these spectra. The maximal ideals of $R_q[G]$
were classified in \cite{Y2}. However, the rest of the inclusions between the prime ideals of $R_q[G]$ are presently unknown and the topology of the spaces 
$\Spec R_q[G]$ remains unknown. 

The quantum group direction of \cite{HL,J}
developed into the study of the spectra of the axiomatic class of algebras $R$ with a rational action of a torus $H$ such that $R$ has finitely 
many $H$-prime ideals. We refer the reader to the book \cite{BG} for a detailed account. Goodearl, Letzter and Brown \cite{BG,GL} obtained a
general stratification result for such algebras which states that $\Spec R$ can be stratified into finitely many strata each of which is homeomorphic 
to a torus (see \S \ref{1.2} for details). The strata are indexed by the $H$-prime ideals of $R$. They are well understood for all quantum groups $R_q[G]$, 
\cite{HLT,J0,Y2} and quantum Schubert cell algebras \cite{BCL,LY,MC,Y1,Y1b}. Finally, in the general setting, Goodearl and Brown \cite{BG2} stated a conjecture that 
would fully describe the topology of the spaces $\Spec R$, i.e., how the above strata of $\Spec R$ are glued together. The conjecture is reviewed in \S \ref{1.3}
and the results in the paper are stated in \S \ref{1.4}.  
\subsection{Stratifications of prime spectra}
\label{1.2}
Throughout, $\KK$ will denote an infinite field of arbitrary characteristic.

Let $R$ be a $\KK$-algebra with a rational action of the $\KK$-torus $H$ by algebra automorphisms.  Call an ideal $H$-prime if it is prime in the set of $H$-invariant ideals of $R$, and denote the set of $H$-prime ideals of $R$ by $\HSpec R$. (The reader is referred to \cite[II.1-2]{BG} for background on $H$-primes and related concepts.)  We restrict our attention to the setting where

(C) the $H$-spectrum of $R$ is finite and $R$ satisfies the noncommutative Nullstellensatz \cite[9.1.4]{McRo}.

Following \cite{BG, GL}, for $I \in \HSpec R$, define the stratum
\[
\Spec_I R := \{ P \in \Spec R \mid \cap_{h \in H} h. P = I \}. 
\] 

\bth{strat} {\em{\cite[Theorem II.2.13]{BG}, \cite[Lemma 3.3]{BG2}}} In the above setting the following hold:

(i) For all $I \in \HSpec R$, the set $\ol{E}_I$ of all regular homogeneous elements in $R/I$ is an Ore set,
and the localization $R_I:= (R/I)[\ol{E}_I^{\; -1}]$ is an $H$-simple ring (i.e., has no non-trivial $H$-primes).

(ii) We have the homeomorphisms $\Spec_I R \cong \Spec R_I \cong \Spec Z(R_I)$ obtained by
localization, contraction, and extension. The center $Z(R_I)$ is a Laurent polynomial ring over $\KK$ 
in at most $\rank H$ generators.

(iii) The result in (ii) is also valid when $\ol{E}_I$ is replaced by any Ore set $E_I \subset \ol{E}_I$ such that  
$(R/I)[E_I^{-1}]$ is $H$-simple. Furthermore, $Z((R/I)[E_I^{-1}]) = Z(R_I)$.

(iv) The decomposition 
\[
\Spec R = \bigsqcup_{I \in \HSpec R} \Spec_I R
\]
forms a finite stratification of $\Spec R$. 
\eth
\subsection{The Brown--Goodearl conjecture}
\label{1.3}
The description of the Zariski topology of $\Spec R$ is equivalent to describing how the different strata in \thref{strat} (iv) are glued together.
Let $CL(T)$ denote the set of closed sets in a topological space $T$.  If $X$ is a closed set in $\Spec_IR$ (with respect to the induced topology from $\Spec R$), write $\overline{X}$ for its closure in $\Spec R$.  For each pair of incident $H$-primes $I \subseteq J$, define a map on the closed sets by
\begin{equation}\label{eq:maps varphi}\varphi_{IJ}: CL(\Spec_IR) \longrightarrow CL(\Spec_JR): X \mapsto \overline{X} \cap \Spec_JR.
\end{equation}
In general, if $Y$ is a closed set in $\Spec R$, then $\varphi_{IJ}(Y \cap \Spec_IR) \subseteq Y \cap \Spec_JR$.  However, if we fix $I \in \HSpec(R)$ and take $P \in \Spec_IR$ and $Y = V(P)$, then it follows easily from the properties of the Zariski topology that this inclusion becomes an equality for all $\varphi_{IJ}$ with $J \supseteq I$.  In particular, we have
\begin{align*}
V(P) & = \{Q \in \Spec R: Q \supseteq P\} \\
&= \bigsqcup_{J \supseteq I} \big(V(P) \cap \Spec_JR\big) \\
& = \bigsqcup_{J \supseteq I} \varphi_{IJ}\big(V(P) \cap \Spec_IR\big).
\end{align*}
Describing $\varphi_{IJ}$ is equivalent to the description of the topological structure of $\Spec R$.  The Brown--Goodearl conjecture aims to give an alternative description of the maps $\varphi_{IJ}$ in terms of certain commutative algebras and homomorphisms between them, with a view to allowing the sets $V(P)$ to be computed explicitly.  In order to state the conjecture, we first need to introduce some notation.

Let $I \subseteq J$ be $H$-primes in $R$, and let $E$ be an Ore set in $R/I$. 
We say that $E$ is a {\em{separating Ore set for the incident pair}} $(I, J)$ if the following two conditions hold:
\begin{enumerate}
\item The image of $E$ in $R/J$ consists of regular elements.
\item If $K$ is an $H$-prime of $R$ such that $I \subseteq K$ and $K\not\subseteq J$, then $E \cap K \neq \varnothing$.
\end{enumerate}  

Take $E_I$ to be any Ore set of regular $H$-eigenvectors in $R/I$ such that the localization $R_I:= (R/I)[E_I^{-1}]$ is $H$-simple (as in \thref{strat} (iii)) 
and define
\begin{equation}\label{eq:def of EIJ and ZIJ}\begin{aligned}
\EE_{IJ} & = \{c \in R/I: c \text{ is an $H$-eigenvector, } c \not\in J/I\},\\
Z_{IJ} & = \{z \in Z(R_I) : zc \in R/I \text{ for some }c \in \EE_{IJ}\}.
\end{aligned}\end{equation}
Brown and Goodearl proved in \cite[Definition~3.8]{BG2} that $Z_{IJ}$ a well-defined $\KK$-algebra, and that it is a subalgebra of $Z(R_I)$.  Further, if $E_{IJ}$ is a separating Ore set for the pair $(I,J)$, then $Z_{IJ}$ becomes more explicit -- it is the center of a localization:
\[
Z_{IJ} = Z((R/I)[E_{IJ}^{-1}]),
\]
\cite[Lemma~3.9]{BG2}.

Intuitively, by studying $Z_{IJ}$ we are restricting our attention to the sub-poset of $H$-$\Spec R$ with minimal element $I$ and maximal element $J$, and hence we can aim to use it to pass from $\Spec_IR$ to $\Spec_JR$.  Let
\begin{equation}\label{eq: def of gIJ}
g_{IJ}: Z_{IJ} \xhookrightarrow{\ \ \ } Z(R_I)
\end{equation}
denote the inclusion map, and let $\pi_{IJ}: R/I \rightarrow R/J$ be the natural projection map.  Define 
\begin{equation}\label{eq: def of fIJ}
f_{IJ} : Z_{IJ}  \longrightarrow  Z(R_J):  z  \mapsto \pi_{IJ}(zc)\pi_{IJ}(c)^{-1},
\end{equation}
where for each $z \in Z_{IJ}$, $c$ is any element of $\EE_{IJ}$ such that $zc \in R/I$.  This is a homomorphism of $\KK$-algebras by \cite[Lemma~3.10]{BG2}; 
when there is a separating Ore set  for the pair $I \subseteq J$ (i.e., $Z_{IJ}$ can be realized as a localization), 
the definition of $f_{IJ}$ simplifies to the natural projection map induced by universality.

Finally, let $g_{IJ}^{\circ}: \Spec Z(R_I) \rightarrow \Spec Z_{IJ}$, $f_{IJ}^{\circ}: \Spec Z(R_J) \rightarrow \Spec Z_{IJ}$ denote the comorphisms induced from $g_{IJ}$, $f_{IJ}$.  We are now in a position to state the Brown-Goodearl conjecture.

\bcj{BGconj} \cite[Conjecture 3.11]{BG2} Assume that $R$ is a $\KK$-algebra with a rational action of the $\KK$-torus $H$ by algebra automorphisms satisfying the condition (C). Let $I \subset J$ be $H$-primes in $R$, and identify $\Spec_IR$, $\Spec_JR$ with $\Spec Z(R_I)$, $\Spec Z(R_J)$.  
Define a map $CL(\Spec_IR) \longrightarrow CL(\Spec_JR)$ by
\[f_{IJ}^{\circ}\overline{|}g_{IJ}^{\circ}(X) = (f_{IJ}^{\circ})^{-1}(\overline{g_{IJ}^{\circ}(X)}),\]
where $\overline{(\cdot)}$ denotes closure in $\Spec Z_{IJ}$.  It is conjectured that
\[
\varphi_{IJ} = f_{IJ}^{\circ}\overline{|}g_{IJ}^{\circ},
\]
where $\varphi_{IJ}$ is the map defined in \eqref{eq:maps varphi}.
\ecj

The conjecture has so far only been verified for a few small algebras, e.g. $R_q[M_2]$ and $R_q[{\mathrm{SL}}_3]$ in \cite{BG2}.  This was later used to study the topology of $\Spec R_q[\mathrm{SL}_3]$ in \cite{F}.  The existence of separating Ore sets has been proved for $R_q[M_{m,n}]$ and $R_q[SL_n]$ in \cite{F2}.
\subsection{Separating Ore sets for all quantum groups and quantum Schubert cell algebras}
\label{1.4} In this paper we construct explicit separating Ore sets for all incident pairs of $H$-prime ideals for the two largest 
classes of quantum algebras that fit the setting of \cjref{BGconj}. The first one is the family of {\em{quantized coordinate rings}} 
of all complex simple Lie groups $R_q[G]$. The second one is the class of {\em{quantum Schubert cell algebras}} ({\em{or quantum unipotent cells}}).
For each symmetrizable Kac--Moody algebra $\g$ and a Weyl group element $w \in W$, the quantum 
Schubert cell algebra $\UU^-[w]$ is a deformation of the universal enveloping algebra $\UU(\n_- \cap w(\n_+))$, 
where $\n_\pm$ are the nilradicals of a pair of opposite Borel subalgebras. The algebras were introduced in \cite{DKP,L}.
The algebras of quantum matrices $R_q[M_{m,n}]$ arise as special cases of the algebras $\UU^-[w]$ for $\g = {\mathfrak{sl}}_{m+n}$ 
and a particular choice of $w$ in the symmetric group $S_{m+n}$. Other interesting families of algebras, such as the quantum Euclidean spaces, 
the algebras of quantum symmetric and antisymmetric matrices arise as special cases of $\UU^-[w]$ from finite dimensional Lie algebras $\g$; 
the affine case of the construction also leads to interesting families of algebras \cite{JN}.  We refer the reader to \S \ref{2.2} and \S \ref{3.1} for details on the
definitions of the algebras $R_q[G]$ and $\UU^-[w]$.

In the above two cases the torus $H$ is a maximal torus of a finite dimensional simple Lie group or a Kac--Moody group, respectively. 
The $H$-primes are parametrized by elements of the Weyl group $W$ of $\g$: 
\[
H\text{-}\Spec R_q[G] = \{ I_{u,v} \mid u, v \in W \}, \quad
H\text{-}\Spec \, \UU^-[w] = \{ I_w(u) \mid u \in W,  u \leq w \},
\] 
and the inclusions are given by 
\[
I_{u',v'} \subseteq I_{u, v} \Leftrightarrow u' \geq u, v' \geq v, \quad
I_w(u) \subseteq I_w(u') \Leftrightarrow u \leq u',
\]
see \cite[Theorem 4.4]{HLT} and \cite[Proposition 10.3.5, Lemma 10.3.6]{J} for the results on $R_q[G]$, and \cite[Theorem 1.1]{Y1} for those on $\UU^-[w]$. Denote by $\PP^+$ the set of 
dominant (integral) weights of $\g$. 
For $\la \in \PP^+$ and $u,v,w \in W$, one defines the 
quantum minors 
\[
\De_{u \la, v \la} \in R_q[G]
\]
and the quasi $R$-matrix
\[
\RR^w \in \UU^+[w] \wh{\otimes} \, \UU^-[w]
\]
for an appropriate completion of the tensor product via the grading by the root lattice \cite[\S 4.1.1]{L}, see \eqref{q-min} and \eqref{Rw} for details.
Denote by $\tau$ the anti-automorphism
of $\UU_q(\g)$ defined by
\begin{equation}
\label{tau}
\tau(E_i) = E_i, \; \; \tau(F_i) = F_i, \; \; \mbox{and} \; \; \tau (K_i) = K_i^{-1}
\end{equation}
on the Chevalley generators of $\UU_q(\g)$, see \S \ref{2.1} for details. The next theorem summarizes our construction of 
separating Ore sets. 
\medskip
\\
\noindent
{\bf{Main Theorem.}}
{\em{Let $\KK$ be an infinite field of arbitrary characteristic and $q \in \KK^{*}$ a non-root of unity.

(i) For every connected, simply connected, complex simple algebraic group $G$, there exist separating Ore sets for all incident 
pairs of $H$-primes of $R_q[G]$. More precisely, for all Weyl group elements $u' \geq u, v' \geq  v$, the multiplicative subset 
$E_{u,v}$  generated by the elements $\De_{u \la, \la}$ and $\De_{-v \mu, - \mu}$ for $\la, \mu \in \PP^+$
is a separating Ore set for the pair $I_{u',v'} \subseteq I_{u, v}$.

(ii) For every symmetrizable Kac--Moody algebra $\g$ and $w \in W$, there exist separating Ore sets for all incident 
pairs of $H$-prime ideals of the quantum Schubert cell algebra $\UU^-[w]$. More precisely, for all $u', u \in W$, $u' \leq u \leq w$,
\[
E_w(u) = q^\Zset \{ (\De_{u \la, w \la} \otimes \id) (\tau \otimes \id) \RR^w : \lambda \in \PP^{+}\}
\]
is a separating Ore set for the pair $I_w(u') \subseteq I_w(u)$, where $q^\Zset = \{ q^n \mid n \in \Zset \}$. 
}}
\medskip

In \cite{F2}, separating Ore sets were constructed for all pairs of incident $H$-prime ideals in the algebras of quantum matrices $R_q[M_{m,n}]$; these sets agree with those constructed in the Main Theorem (ii) above (restricted to the appropriate choice of $\g$ and $w$).  Readers interested in realizing the separating Ore sets in ring-theoretic terms for algebras of quantum matrices are referred to \cite[\S4]{F2} for more details and examples.

In \cite{GY} quantum cluster algebra structures were constructed on the quantized coordinate rings of all double Bruhat cells
and in \cite{LY} large families of quantum seeds were constructed on the quantized coordinate rings of open Richardson varieties. 
These constructions and the explicit Ore sets from the Main Theorem can be used to construct cluster algebra models for all
localizations $(R/I)[E_{IJ}^{-1}]$ entering in the Brown--Goodearl conjecture for the algebras $R_q[G]$ and $\UU^-[w]$.
 
For both families of algebras, the constructed separating Ore sets $E_{IJ} \subset R/I$ in the Main Theorem and 
the Ore sets $E_J \subset R/J$, satisfying the condition in \thref{strat} (iii), 
exhibit two important special properties (see \S \ref{3.4} for details): For all $I \subset J \in \HSpec R$,

(1) the sets $E_{IJ}$ and $E_J$ are the projections of Ore sets in
$R$ (to be denoted also by $E_J$ and $E_{IJ}$) and

(2) $E_{IJ} = E_J$.
 
These properties and the above quantum cluster algebra models shed light on the nature of the maps $\varphi_{IJ}$ and $f_{IJ}^{\circ}\overline{|}g_{IJ}^{\circ}$.
We expect that they will play a key role in the resolution of the Brown--Goodearl conjecture.

We finally note that \cite[\S4.5]{Y2} constructed a Dixmier type map 
\[
\Dix_G \colon \Symp (G, \pi) \to \Prim R_q[G]
\]
from the symplectic foliation of the standard Poisson structure on $G$ (see \cite{Dr} and \cite[Theorem A.2.1]{HL}) to the primitive spectrum of $R_q[G]$. It was 
proved in \cite[Theorem 4.6]{Y2} that the map is an $H$-equivariant bijection, and was conjectured that it is a homeomorphism. This would 
crystallize the Hodges--Levasseur--Joseph orbit method program for $R_q[G]$. The construction of the separating Ore sets 
for $R_q[G]$ also presents a step towards the realization of this program. 
\medskip
\\
\noindent
{\bf Acknowledgements.} S. F. would like to thank Robert Marsh for 
many clear and patient explanations on the general theory and background of quantum groups.
M. Y. would like to thank Newcastle University and the Max Planck Institute for Mathematics in Bonn
for the hospitality during visits in the Fall of 2015 when the results in the paper were obtained.  
We are also thankful to the referee whose valuable suggestions helped us to improve the exposition.
\sectionnew{Separating Ore sets for Quantum groups}
\lb{Ore-quantGr}
In this section we construct explicit separating Ore sets for all incident pairs of torus-invariant 
prime ideals of the quantum function algebras on complex simple Lie groups $R_q[G]$, proving part (i) of the 
Main Theorem.
The section also contains results on Kac--Moody algebras needed in the next section.
\subsection{Quantized universal enveloping algebras}
\label{2.1}
We fix a symmetrizable Kac--Moody algebra 
$\g$ of rank $r$ with Cartan matrix $(c_{ij})$ and Weyl group $W$. Let $\UU_q(\g)$ be 
the quantized universal enveloping algebra of $\g$ over an infinite  
field $\KK$ of arbitrary characteristic. The deformation parameter $q$ will be assumed to be a non-root of unity.  Here and below, for $m \leq n \in \Zset$, we set $[m,n] :=\{m, \ldots, n \}$. 
The Chevalley generators of $\UU_q(g)$ will be denoted by 
\[
E_i, F_i , K_i^{\pm 1}, \; \; i \in [1,r].
\]
Denote by $\UU^+_q(\g)$ and $\UU^-_q(\g)$ the subalgebras
of $\UU_q(\g)$ generated by $\{E_i\}$ and $\{F_i\}$ respectively.

We will follow the notation of \cite{Ja}, except that the generators of $\UU_q(\g)$ will be indexed by $[1,r]$ instead of the 
simple roots of $\g$. The root and weight lattices of $\g$ will be denoted by 
$\QQ$ and $\PP$, and the set of dominant integral weights of $\g$ by $\PP^+$. 
Let $\al_1, \ldots, \al_r$ and $\al_1\spcheck, \ldots, \al_r\spcheck$ be the simple roots and coroots  of $\g$, 
and $\lcor.,. \rcor$ be the symmetric bilinear form on $\h= \Span \{ \al_i \}_{i=1}^r$ such that 
$\lcor \al_i, \al_j \rcor = d_i c_{ij}$,
$\forall i, j \in [1,r]$. Denote by $\PP^{++}$ the set of strictly dominant weights of $\g$, 
\[
\PP^{++} = \{ \la \in \PP \mid \lcor \la, \al_i\spcheck \rcor > 0, \; \; \forall i \in [1,r] \}.
\]
Let $\vpi_1, \ldots, \vpi_r$ be the fundamental weights of $\g$.
The weight spaces of a $\UU_q(\g)$-module $V$ 
are defined by 
\[
V_\nu = \{ v \in V \mid K_i v = q^{ \lcor \nu, \al_i \rcor} v, \; 
\forall i \in [1,r] \}, \; \nu \in \PP.
\]
For $\la \in \PP^+$, let $L(\la)$ be the unique irreducible 
$\UU_q(\g)$-module of highest weight $\la$. Let $L(\la)\spcheck$ be its graded
dual which is a $\UU_q(\g)$-module via
\[
\lcor x \xi, b \rcor := \lcor \xi, S(x) b \rcor, \quad x \in \UU_q(\g),\ b \in L(\la),\ \xi \in L(\la)\spcheck, 
\]
where $S$ denotes the antipode of $\UU_q(\g)$.  We will use the compatible actions 
of the braid group $\B_\g$ on $\UU_q(\g)$ and $L(\la)$,  
\cite[\S 8.6 and \S 8.14]{Ja}. As usual, $q$-integers and $q$-factorials will be denoted by
\begin{equation}\label{q-factorials}
[n]_q := \frac{q^n - q^{-n}}{q-q^{-1}}, \quad [n]_q! := [1]_q[2]_q\dots [n]_q.
\end{equation}
\subsection{Quantum function algebras}
\label{2.2} 
Denote by $G$ the complex Kac--Moody group associated to $\g$. In the finite dimensional case, 
this is the connected simply connected algebraic group with Lie algebra $\g$. 

For any infinite field $\KK$ and a non-root of unity $q \in \KK^*$, 
one defines the quantum coordinate ring $R_q[G]$ as the subalgebra of $(\UU_q(\g))^*$ 
spanned by the matrix coefficients of the modules $L(\la)$, $\la \in \PP^+$. They are denoted by
\begin{equation} 
\label{c-notation}
c_{\xi, b} \in (\UU_q(\g))^*,\quad 
c_{\xi, b}(x) = \lcor \xi, x b \rcor, \quad 
\xi \in L(\la)\spcheck, b \in L(\la), 
x \in \UU_q(\g).
\end{equation}
The algebra $R_q[G]$ is a quantum analog of the algebra of strongly regular functions on $G$, \cite[Sect. 4]{KP}.
With respect to the canonical action of $\UU_q(\g) \otimes \UU_q(\g)$ on $R_q[G]$ we have the isomorphism
\[
R_q[G] \cong \bigoplus_{\la \in \PP^+} L(\la)\spcheck \otimes L(\la).
\]
In particular $R_q[G]$ is $\PP \times \PP$ graded by
\[
R_q[G]_{\mu, \nu} = \{ c_{\xi, b} \mid \la \in \PP^+, \xi \in (L(\la)\spcheck)_\mu, b \in L(\la)_\nu \}, \quad \mu, \nu \in \PP.
\]
Identify the (rational) character lattice of $H:= (\KK^*)^r$ with $\PP$ via
\begin{equation}
\label{char-latt}
\nu \in \PP \quad \mbox{maps to the character} \quad 
t=(t_1, \ldots, t_r) \mt  t^\nu:=\prod_k t_k^{ \lcor \nu, \al_i\spcheck \rcor}.
\end{equation}
The torus $H$ acts on $R_q[G]$ on the left and the right by algebra automorphisms by 
\begin{equation}
\label{H-act}
t \cdot c := t^\nu c,   \quad
c \cdot t  := t^\mu c \quad
\mbox{for} \; \; c \in R_q[G]_{\nu, \mu},\ t \in H.
\end{equation}

For each $\la \in \PP^+$, we fix a highest weight vector $b_\la$ of $L(\la)$. For $u \in W$, set $b_{u \la} := 
T_{u^{-1}}^{-1} b_\la$. It is well known that the vector $b_{u \la}$ only depends on $u\la$ and not on the choice of $u$. 
Let $\xi_{u \la} \in (L(\la)\spcheck)_{- u\la}$ 
be the unique vector such that $\lcor \xi_{u \la}, b_{u \la} \rcor =1$. For $u, v \in W$, define the 
quantum minors
\begin{equation}
\label{q-min}
\De_{u \la, v \la} := c_{\xi_{u \la}, b_{v \la}} \in R_q[G].
\end{equation}
One has
\begin{equation}
\label{prod-minor}
\De_{u \la, v \la}  \De_{u \mu, v \mu}  = \De_{u (\la + \mu), v (\la + \mu)}, \quad
\la, \mu \in \PP^+.
\end{equation}
Define the subalgebra
\[
R^+ = \Span \{ c_{\xi, b_\la} \mid \la \in \PP^+, \xi \in L(\la)\spcheck \}
\]
of $R_q[G]$. It is a quantum analog of the coordinate ring of $G/U^+$ where $U^\pm$ are the unipotent radicals of the 
pair of opposite Borel subgroups corresponding to the fixed triangular decomposition of $\g$ used in defining $\UU_q(\g)$. 
Define the multiplicative subsets
\[
E_u^+ := \{  \De_{u \la, \la} \mid \la \in \PP^+ \} \subset R^+
\]
(cf. \eqref{prod-minor}) and the ideals
\[
I_u^+ = \Span \{ c_{\xi, b_\la} \mid \la \in \PP^+, \xi \in L(\la)\spcheck, \,
\xi \perp \UU^+_q(\g) b_{u \la} \} \subseteq R^+.
\]
\ble{aux} For all symmetrizable Kac--Moody algebras $\g$ and $u \in W$,

(i) $E_u^+$ is an Ore set in $R^+$,

(ii) $I_u^+$ is a completely prime ideal of $R^+$, $E_u^+ \cap I_u^+ = \varnothing$, 
and the image of $E_u^+$ in $R^+/I_u^+$ consists of normal elements:
\[
\De_{u \la, \la} c = q^{\lcor u \la, \mu \rcor - \lcor \la, \nu \rcor} c \De_{u \la, \la} \mod I_u^+, \quad 
\forall c \in R^+_{- \mu, \nu},\  \mu \in \PP,\ \la, \nu \in \PP^+.
\]
\ele
The second part of the lemma was stated and proved in \cite{HL,J} for finite dimensional groups
$G$. However those proofs and the one in \cite[\S 3.2]{Y1b} easily extend to the Kac--Moody case. 
The proof of the first part in \cite[Lemma 9.1.10]{J} used an induction from the longest element of $W$ down 
to $u \in W$, and in that form it applies to finite dimensional groups $G$. In \S\ref{skew derivations} below we give a proof by a reverse 
induction from $1 \in W$ to $u \in W$ that works for all symmetrizable Kac--Moody algebras $\g$.

For the remaining part of this section we restrict to finite dimensional simple Lie algebras $\g$. Consider
the lowest weight vector $b_{-\la}:= T_{w_0} b_{- w_0 \la}$ of  $L_{ - w_0 \la}$ (where $w_0$ is the longest element of $W$) and define the subalgebra 
\[
R^- = \Span \{ c_{\xi, b_{- \la}} \mid \la \in \PP^+ \} 
\]
of $R_q[G]$. It is a quantum analog of the coordinate ring of $G/U^-$. For $u, v \in W$, define the ideals
\[
I_v^- := \Span \{ c_{\xi, b_{- \la} } \mid \la \in \PP^+, \xi \in L(\la)\spcheck, \,
\xi \perp \UU^-_q(\g) b_{- v \la} \}\ \subseteq \ R^-
\]
and 
\[
I_{u,v} := I_u^+ R^- + R^+ I_v^- \subseteq R_q[G].
\]
Denote the multiplicative subset
\[
E_v^- := \{  \De_{- v \la, - \la} \mid \la \in \PP^+ \} \subset R^-,
\] 
recall \eqref{prod-minor}. Note that the quantum minors $\De_{-v \la, - \la}$ come from the matrix coefficients 
of the irreducible highest weight modules $L(- w_0 \la)$.

Denote by $E_{u,v}$ the multiplicative subset of $R_q[G]$ generated by $E_u^+$ and $E_v^-$.
\bre{fund1} Because of \eqref{prod-minor}, $E_{u,v}$ can be equivalently defined as the multiplicative subset of $R_q[G]$ generated by 
the set of quantum minors coming from fundamental representations
\[
\De_{u \vpi_i, \vpi_i},  \De_{- v \vpi_i, - \vpi_i}, \quad i \in [1,r].
\]
\ere
\subsection{Skew derivations and Ore sets}\label{skew derivations}
In order to prove \leref{aux} (i) we first need to introduce some new notation.  Let $\sig$ be an automorphism of a $\KK$-algebra $A$ and 
$\del$ be a locally nilpotent right or left $\sig$-derivation 
of $A$ such that $\sigma \del \sigma^{-1} = q' \del$ for some $q' \in \KK^*$.
The degree $\deg_\del a$ of an element $a \in A \backslash \{ 0 \}$ 
is defined to be the minimal positive integer $m$ such that 
$\del^{m+1}(a) =0$. 
For two $\sig$-eigenvectors $a, b \in A$ of $\del$-degrees $m$ and $n$ respectively, write
\[
\del^k (a b) = \sum_{i=0}^{m+n-k} s_{k,i} ( \del^{m-i} a) ( \del^{k-m+i} b) \quad \mbox{with} \quad
s_0, \ldots, s_{m+n-k} \in \KK.
\]
The skew derivation $\del$ is called regular if $s_{k,0} \neq 0$ for $k \in [m,m+n]$
and $s_{k,m+n-k} \neq 0$ for $k \in [n,m+n]$. One of the consequences of this is 
that $\deg_\del(ab) = \deg_\del(a) + \deg_\del(b)$. 
If the scalar $q'$ is a non-root of unity, then $\del$ is regular
because all coefficients $s_{k,i}$ are products 
of $q'$-binomial coefficients and $\sig$-eigenvalues.  
\medskip
\\
\noindent
{\em{Proof of \leref{aux} (i).}} The statement holds for the set $E_1$ since it consists of normal elements of $R^+$,
\[
\De_{\la, \la} c = q^{\lcor \la, \nu - \mu \rcor} c \De_{\la, \la}, \quad 
\forall c \in R^+_{- \mu, \nu},\ \mu \in \PP,\ \la, \nu \in \PP^+.
\]
This follows from the $\RR$-matrix commutation relations in $R_q[G]$, see e.g. \cite[Lemma 2.2 (i)]{Y2}. 
By induction on the length of $u$ we prove that,
  
(*) {\em{For $u \in W$ and $\la \in \PP^+$, 
$\{ \De_{u \la, \la}^n \mid n \in \Nset \}$ is an Ore set in $R^+$.}} 
\\
We use the following result of Joseph \cite[Lemma A.2.9]{J}:

(**) {\em{For a regular skew derivation $\del$ of an algebra $A$, if $\{ e^n \mid n \in \Nset \}$ 
is an Ore set in $A$ and $e \in A$ is a $\sig$-eigenvector of degree $\deg_\del e = m$, 
then $\{ (\del^m e)^n \mid n \in \Nset \}$ is an Ore set in $A$}}. 

$R^+$ is a $\UU_q(\g)$-module algebra by 
\[
x \cdot c_{\xi, b_\la} := c_{x \cdot \xi, b_\la}, \quad x \in \UU_q(\g), \xi \in L(\la)\spcheck, \la \in \PP^+.
\]
We apply (**) for $\sig_i := (K_i \cdot)$ and $\del_i := (E_i \cdot)$. The skew derivations $\del_i$ are 
locally nilpotent because $L(\la)$ is integrable. They are regular since $\sig_i \del_i \sig_i^{-1} = q_i^2 \del_i$. 

Let $u \in W$ and $i \in [1,r]$ be such that $\ell(s_i u) = \ell(u) +1$ where $\ell \colon W \to \Nset$ is the 
length function. Then $b_{u \la}$ is a highest weight vector for the $\UU_{q_i} (\sl_2)$ subalgebra 
of $\UU_q(\g)$ generated by $E_i, F_i$ and $K_i$ of highest weight $\lcor u \la, \al_i\spcheck \rcor \vpi_i$ 
and
\[
b_{s_i u \la} = T_i^{-1} b_{u \la} = F_i^{\lcor u \la, \al_i\spcheck \rcor} b_{u \la}/[\lcor u \la, \al_i\spcheck \rcor]_{q_i}!. 
\]
This implies that 
\[
\deg_{\del_i} \De_{u \la, \la} = \lcor u \la, \al_i\spcheck \rcor \quad \mbox{and} \quad
\del_i^{\lcor u \la, \al_i\spcheck \rcor} \De_{u \la, \la} = t \De_{s_i u \la, \la}
\]
for some $t \in \KK^*$. Now (*) follows from this and (**) by induction on $\ell(u)$.
\qed 
\medskip
\\

\subsection{Stratifications of $\Spec R_q[G]$}
\label{2.3}

The following theorem was proved by Hodges--Levasseur--Toro \cite{HL,HLT} and Joseph \cite{J0,J}, see
\cite[Proposition 8.9]{J0}, \cite[Proposition 10.1.8, Proposition 10.3.5]{J}. 

\bth{H-primes} \cite{HLT,J0} Let $\g$ be a complex simple Lie algebra, $\KK$ an arbitrary infinite field and $q \in \KK^*$ a non-root of unity.

(i) For each pair $(u, v) \in W \times W$, $I_{u,v}$ is an $H$-invariant completely prime ideal of $R_q[G]$ with respect to the left action of $H$. 
All $H$-primes of $R_q[G]$ are of this form. The ideals $I_{u,v}$ are also invariant under the right $H$-action and all $H$-primes of $R_q[G]$ 
with respect to the right action are of this form. For all $u, v, u', v' \in W$, 
\[
I_{u', v'} \subseteq I_{u, v} \quad  \Leftrightarrow \quad u' \geq u \; \; \mbox{and} \; \; v' \geq v.
\]

(ii) For $(u, v) \in W \times W$, $E_{u,v}$ is an Ore set in $R_q[G]$, $E_{u,v} \cap I_{u,v} = \varnothing$ and the image of $E_{u,v}$ in $R_q[G]/I_{u,v}$ consists of normal elements,
\begin{align}
\label{n1}
\De_{u \la, \la} c &= q^{\lcor u \la, \mu \rcor - \lcor \la, \nu \rcor} 
c \De_{u \la, \la} \mod I_{u,v}, 
\\
\label{n2}
\De_{-v \la, - \la} c &= q^{ \lcor v \la, \mu \rcor - \lcor \la, \nu \rcor} 
c \De_{-v \la, - \la} \mod I_{u,v} 
\end{align}
for all $\la \in \PP^+$, $\mu, \nu \in \PP$ and $c \in R_q[G]_{-\mu, \nu}$.
The localization $(R_q[G]/I_{u,v}) [E_{u,v}^{-1}]$ is $H$-simple. 
\eth

In \cite{HLT,J0}, \thref{H-primes} was stated for $\KK = \Cset$ and a non-root of unity $q \in \Cset^*$ and in \cite{J} it was stated for $\KK = k(q)$ 
where $k$ is a field of characteristic 0. It is well known that the proofs in \cite{HLT,J0,J} work in the generality stated in \thref{H-primes}.

From now on, when referring  to $H$-invariance of prime ideals of $R_q[G]$, we will have in mind the left action \eqref{H-act} 
of $H$, keeping in mind that the right invariance leads to the same set of ideals.

\bre{prod} It follows from \eqref{prod-minor}, \eqref{n1}, and \eqref{n2} that the image of the Ore set $E_{u,v}$ in $R_q[G]/I_{u,v}$ has the form
\begin{equation}
\label{wrong-set}
\{ q^n \De_{u \la, \la} \De_{-v \mu, - \mu} \mid \la, \mu \in \PP^+, n \in \Sig \}
\end{equation}
for some additively closed subset $\Sig$ of $\Zset$. However, for general $u' \geq u, v' \geq v$, 
the image $E_{u,v}$ in $R_q[G]/I_{u',v'}$ does not have this factorization property. In other words, the subset 
\eqref{wrong-set} of $R_q[G]/I_{u',v'}$ is not a multiplicative subset. 
\ere

Joseph and Hodges--Levasseur--Toro proved the following decomposition theorem for the spectra of 
$R_q[G]$, see \cite[Theorems 8.11, 9.2]{J0} and \cite[Theorems 4.4, 4.15, 4.16]{HLT}.
It also follows from Theorems \ref{tstrat} and \ref{tH-primes}.

\bth{Spec} \cite{HLT,J0} In the setting of \thref{H-primes}, we have the decomposition,
\[
\Spec R_q[G] = \bigsqcup_{u, v \in W} \Spec_{u,v} R_q[G]
\]
where
\[
\Spec_{u,v} R_q[G] = \{ P \in \Spec R_q[G] \mid \cap_{t \in H}\ t \cdot P = I_{u,v} \}.
\]
For all $u, v \in W$, localization, contraction, and extension induce the homeomorphisms 
\[
\Spec_{u,v} R_q[G] \cong \Spec (R_q[G]/I_{u,v}) [E_{u,v}^{-1}] \cong
\Spec Z( (R_q[G]/I_{u,v}) [E_{u,v}^{-1}]).
\]
\eth
\subsection{Separating Ore sets for quantum groups}
\label{2.4} 
The next theorem constructs separating Ore sets for all incident pairs of $H$-prime ideals of $R_q[G]$.  This proves part (i) of the Main Theorem.

\bth{separateOre} Assume that $\g$ is a complex simple Lie algebra, $\KK$ an arbitrary infinite field, $q \in \KK^*$
a non-root of unity. For all $u,v,u',v' \in W$ such that $u' \geq u$ and $v' \geq v$, 
the set $E_{u,v}$ is a separating Ore set for the pair of $H$-prime ideals $I_{u',v'} \subseteq I_{u,v}$ of $R_q[G]$.
\eth

The fact that the above pairs of ideals exhaust all incident pairs $I \subseteq J$ of $H$-primes of $R_q[G]$
follows from \thref{H-primes} (i).
\begin{proof} By \thref{H-primes} (ii), $E_{u,v} \cap I_{u,v} = \varnothing$ and by \thref{H-primes} (i), 
$I_{u',v'} \subset I_{u,v}$. Thus $E_{u,v} \cap I_{u',v'} = \varnothing$. Since $I_{u',v'}$ is a completely prime ideal,
the image of the set $E_{u,v}$ in $R_q[G]/I_{u',v'}$ consists of regular elements. By \thref{H-primes} (ii), $E_{u,v}$ is an Ore set in
$R_q[G]$. Thus, its image in $R_q[G]/I_{u',v'}$ is an Ore set in $R_q[G]/I_{u',v'}$ that does not intersect $I_{u,v} /I_{u',v'}$. 
This establishes the first property for a separating Ore set (cf. \S \ref{1.3}).

Now let $I_{u',v'} \subseteq I_{u'', v''}$ be another incident pair such that $I_{u'',v''} \not\subseteq I_{u,v}$. By \thref{H-primes} (iii), either
$u'' \not\geq u$ or $v''\not\geq v$. Assume that the first condition holds; the second case is treated similarly. 
Let $\la \in \PP^{++}$ be a strictly dominant weight. The inclusion relations between Demazure modules \cite[Proposition 4.4.5]{J} imply that 
\[
L(\la)_{u \la} \not\subseteq \UU_q^+(\g) L(\la)_{u'' \la}. 
\]
Therefore $\De_{u \la, \la} \in E_{u,v} \cap I_{u'',v''}$ which establishes the second property for a separating Ore set (see \S \ref{1.3}).
\end{proof}

\bre{qcluster-RqG} For all pairs $u, v \in W$, the localized algebra $(R_q[G]/I_{u,v})[E_{u,v}^{-1}]$ from \thref{H-primes} (iv) 
is isomorphic to the 
quantized coordinate ring of the double Bruhat cell $G^{u,v}:=B^+ u B^+ \cap B^- v B^-$. It was proved in \cite[Main Theorem]{GY} 
that $(R_q[G]/I_{u,v})[E_{u,v}^{-1}]$ has a canonical structure of a quantum cluster algebra for which all frozen cluster variables are inverted.
These are the algebras that govern the stratification picture in \thref{strat}.

On the other hand, for $u' \geq u \in W$, $v' \geq v \in W$, the separating Ore sets from \thref{separateOre} give rise to the partial 
localizations 
\begin{equation}
\label{part-local}
(R_q[G]/I_{u',v'})[E_{u,v}^{-1}].
\end{equation}
These are the localizations that govern the maps in the Brown--Goodearl Conjecture. 
They are obtained from the quantum cluster algebra structures on $R_q[G^{u',v'}]$ from \cite{GY} by uninverting the frozen cluster variables 
and then localizing at $E_{u,v}$. The latter set is not generated by cluster variables in the initial seed in $R_q[G^{u',v'}]$. However, these are 
cluster algebras of infinite type (except for very special cases). We expect that the generators $\De_{u \vpi_i, \vpi_i}$ and $\De_{- v \vpi_i, - \vpi_i}$ 
of $E_{u,v}$ (recall \reref{fund1}) are cluster variables of $R_q[G^{u',v'}]$ of non-initial seeds. This can be verified 
directly when $R_q[G^{u',v'}]$ is of finite type. If this fact is established, then the partial localizations \eqref{part-local} will 
be realized explicitly as localizations of quantum cluster algebras by unfrozen cluster variables and the maps
$\varphi_{IJ}$ and $f_{IJ}^{\circ}\overline{|}g_{IJ}^{\circ}$ in the Brown--Goodearl conjecture will be realized explicitly 
in terms of cluster algebras.
\ere
\sectionnew{Separating Ore sets for the quantum Schubert cell algebras}
\lb{Ore-Uw}
In this section we construct explicit separating Ore sets for all incident pairs of torus-invariant 
prime ideals of the quantum Schubert cell algebras associated to all symmetrizable 
Kac--Moody algebras.
\subsection{Quantum Schubert cell algebras}
\label{3.1}
Throughout the section $\g$ will denote a symmetrizable Kac--Moody algebra, $\KK$ an infinite base field, and $q$ 
a non-root of unity. Fix a Weyl group element $w \in W$. As in \S\ref{2.2}, we set $H := (\KK^*)^r$.

For a reduced expression $w = s_{i_1} \ldots s_{i_l}$, the roots of the Lie algebra $\n_+ \cap w(\n_-)$ are 
$\{ \be_k := s_{i_1} \ldots s_{i_{k-1}} (\al_k) \mid k \in [1,l]\}$. Denote the root vectors 
\[
E_{\be_k} :=
T_{i_1} \ldots T_{i_{k-1}} 
(E_{i_k}) \quad \mbox{and} \quad 
F_{\be_k} := 
T_{i_1} \ldots T_{i_{k-1}} 
(F_{i_k}).
\]
The quantum Schubert cell algebras $\UU^\pm[w]$ are the subalgebras of 
$\UU^\pm_q(\g)$ generated by $\{E_{\be_1}, \ldots, E_{\be_l} \}$ and 
$\{F_{\be_1}, \ldots, F_{\be_l} \}$, respectively. They are independent of the choice of reduced expression of $w$, 
\cite[Proposition 4.2.1]{L}. The algebra $\UU_q(\g)$ is $\QQ$-graded by $\deg E_i = - \deg F_i = \al_i$, $\deg K_i =0$. The corresponding 
graded components will be denoted by $\UU_q(\g)_\nu$, $\nu \in \QQ$.
Define the $H$-action on $\UU_q(\g)$ 
\[
t \cdot x = t^\nu x, \quad t \in H,\ x \in \UU_q(\g)_\nu,\ \nu \in \QQ 
\]
in terms of the identification \eqref{char-latt} of the character lattice of $H$ with $\PP$. The subalgebras $\UU^\pm[w]$ are stable under this action.

Recall \leref{aux} and define the Joseph algebras \cite[\S 10.3]{J}
\[
S^+_w : = \big( (R^+/I_w^+)[(E_w^+)^{-1}] \big)^H
\]
where the invariant part is computed with respect to the right action \eqref{H-act}. The torus $H$ acts on $S^+_w$ by algebra automorphisms
via the left action \eqref{H-act}. For $u \in W$, $u \leq w$, define the ideals
\[
Q^-_w(u) := 
\{ (c_{\xi, b_\la} + I_w^+) \De_{w \la, \la}^{-1} \mid \la \in \PP^+, \xi \in L(\la)\spcheck, \,
\xi \perp \UU^-_q(\g) b_{u \la} \}\ \subset \ S^w_+.
\]
\subsection{Spectra of the quantum Schubert cell algebras}
\label{3.2}
The $H$-primes of $\Spec \, \UU^-[w]$ are classified using Demazure modules as follows.
Denote $W^{\leq w} := \{u \in W \mid u \leq w \}$, and recall the definition of the $q$-factorial $[n]_q!$ from \eqref{q-factorials}. The quasi $\RR$-matrix for $w$ is given by
\begin{equation}
\label{Rw}
\RR^w := \sum_{m_1, \ldots, m_l \, \in \, \Nset}
\left( \prod_{j=1}^l
\frac{ (q_{i_j}^{-1} - q_{i_j})^{m_j}}
{q_{i_j}^{m_j (m_j-1)/2} [m_j]_{q_{i_j}}! } \right) 
E_{\be_l}^{m_l} \ldots E_{\be_1}^{m_1} \otimes  
F_{\be_l}^{m_l} \ldots F_{\be_1}^{m_1}.
\end{equation}
It lies in the completion of $\UU^+[w] \otimes \UU^-[w]$, \cite[\S 4.1.1]{L}.

Recall the anti-automorphism $\tau$ of $\UU_q(\g)$ given by \eqref{tau}.

\bth{Class} \cite{Y1} Let $\g$ be a symmetrizable Kac--Moody algebra, $w\in W$ a Weyl group element, $\KK$ an infinite base field and $q \in \KK^*$ 
a non-root of unity.

(i) The map
\[
\phi_w \colon 
S^+_w \to \UU^-[w], \quad
\phi_w \big( (c_{\xi, b_\la} + I_w^+) \De_{w \la, \la}^{-1} \big) 
:= \big( c_{\xi, b_{w \la}} \otimes \id \big) (\tau \otimes \id) 
\RR^w, \; \; 
\la \in \PP^+, \xi \in L(\la)\spcheck
\]
is a well defined $H$-equivariant algebra anti-isomorphism, where on the first algebra $H$ acts via the left action \eqref{H-act}.

(ii) For $u \in W^{\leq w}$, 
\[
I_w(u) = \vp_w(Q_w^-(u))
\]
are distinct, $H$-invariant, completely prime  
ideals of $\UU^-[w]$. All $H$-prime ideals of $\UU^-[w]$
are of this form.

(iii) The map $u \in W^{\leq w} \mt I_w(u) \in \HSpec \UU^-[w]$ 
is an isomorphism of posets with respect to the Bruhat order 
on $W^{\leq w}$ and the inclusion order on $\HSpec \UU^-[w]$.

(iv) For $\la \in \PP^+$, define
\[d_{u, \la} := \phi_w( (\De_{u \la, \la} +I_w^+) \De_{w \la, \la}^{-1}).\]
The set 
\[
E_w(u) := q^\Zset \{ d_{u, \la}  \mid \la \in \PP^+ \}
\]
is an Ore set in $\UU^-[w]$, and $E_w(u) \cap I_w(u) = \varnothing$. The image of $E_w(u)$ in $\UU^-[w]/I_w(u)$ 
consists of normal elements:
\[
d_{u,\la} c = q^{ \lcor (w+u) \la, \nu \rcor } c d_{u,\la}, \quad \la \in \PP^+,\ c \in \UU^-[w]_\nu,\ \nu \in \QQ,
\]
and the localization $(\UU^-[w]/I_w(u))[E_w(u)^{-1}]$ is $H$-simple.
\eth
In \cite{Y1} the theorem was stated for complex simple Lie algebras and for fields of characteristic $0$. 
All steps in the proof, together with the needed results from \cite{G,J}, appear 
in \cite[Sect. 2--3]{Y1b}, \cite[\S 2.8]{Y2} where it was shown that the proof extends to base fields $\KK$ of arbitrary characteristic. 
The proofs in \cite[Sect. 2--3]{Y1b}, \cite[\S 2.8]{Y2} apply in the more general case of symmetrizable Kac--Moody algebras $\g$, with the
only needed modification provided by \leref{aux} (i). 

\bre{fund2} Eq. \eqref{prod-minor}, \leref{aux} (ii), and \thref{Class} (i) imply that
\[
d_{u, \la}  d_{u, \mu} = q^n d_{u, \la+\mu}, \quad
\la, \mu \in \PP^+
\]
for some $n \in \Zset$ depending on $\la$ and $\mu$. Firstly, this is why $E_w(u)$ is a multiplicative subset of $\UU^-[w]$. 
Secondly, it follows from this identity that $E_w(u)$ can be equivalently 
defined as the multiplicative subset of $\UU^-[w]$ generated by $d_{u, \vpi_i}$ for $i \in [1,r]$ and 
the scalars $q^{\pm1}$, the point being that the generators of $E_w(u)$ come from the fundamental representations of $\UU_q(\g)$. 
\ere

Combining \thref{strat} and \thref{Class}, we obtain:
\bco{stratUw} In the setting of \thref{Class} we have the decomposition
\[
\Spec \, \UU^-[w] = \bigsqcup_{u \in W^{\leq w}} \Spec_u \, \UU^-[w]
\]
where
\[
\Spec_u \, \UU^-[w] = \{ P \in \Spec \, \UU^-[w] \mid \cap_{t \in H}\ t \cdot P = I_w(u) \}.
\]
For all $u \in W^{\leq w}$, localization, contraction, and extension induce the homeomorphisms 
\[
\Spec_u \, \UU^-[w] \cong \Spec (\UU^-[w]/I_w(u)) [E_w(u)^{-1}] \cong
\Spec Z( (\UU^-[w]/I_w(u)) [E_w(u)^{-1}] ).
\]
\eco
The dimensions of the strata were determined in \cite{BCL,Y1b}.
\subsection{Separating Ore sets for the quantum Schubert cell algebras}
\label{3.3}
We are now in a position to prove part (ii) of the Main Theorem.
\bth{SetOreUw}  Assume that $\g$ is a symmetrizable Kac--Moody algebra, $w \in W$ a Weyl group element, $\KK$ is an infinite base field and $q \in \KK^{*}$ 
a non-root of unity.

(i) For all $u \in W^{\leq w}$, the set $E_w(u)$ is an Ore set in $\UU^-[w]$. 

(ii) For $u' \leq u \in W^{\leq w}$, $E_w(u)$ is a separating Ore set for the incident pair $I_w(u') \subseteq I_w(u)$ 
of $H$-primes of $\UU^-[w]$. 
\eth
\thref{Class} (ii)-(iii) implies that the pairs in part (ii) exhaust all incident pairs of $H$-primes of $\UU^-[w]$.  
\begin{proof} (i) By \leref{aux}, $E_u^+$ is an Ore set in $R^+$ and $E_u^+ \cap I_w^+ = \varnothing$ 
because of the inclusion of Demazure modules \cite[Lemma 4.4.3 (v)]{J},
\[
L(\la)_{u \la} \subseteq \UU_q^+(\g) L(\la)_{w \la}, \quad \la \in \PP
\]
(using that $u \leq w$). It follows from the commutation relation in \leref{aux} that 
\[
q^\Zset \{ (\De_{u \la, \la} +I_w^+) \De_{w \la, \la}^{-1} \mid \la \in \PP^+ \}
\]
is an Ore set in $S^+_w$. Part (i) now follows from the anti-isomorphism in \thref{Class} (i).

(ii) By \thref{Class} (iv), $E_w(u) \cap I_w(u) = \varnothing$ and by \thref{Class} (iii), 
$I_w(u') \subseteq I_w(u)$. Therefore, $E_w(u) \cap I_w(u') = \varnothing$. Since the ideal $I_w(u')$ is completely prime,
the image of the set $E_w(u)$ in $\UU^-[w]/I_w(u')$ consists of regular elements. The first part of the theorem implies that 
the image $E_w(u)$ is an Ore set in $\UU^-[w]/I_w(u')$ that does not intersect $I_w(u)/I_w(u')$. 
This proves that $E_w(u)$ satisfies the first property for a separating Ore set.

Fix another incident pair of $H$-primes of $\UU^-[w]$: $I_w(u') \subseteq I_w(u'')$ such that $I_w(u'') \not\subseteq I_w(u)$. It follows from \thref{Class} (iii)
that $u'' \not\leq u$. Consider a strictly dominant weight $\la \in \PP^{++}$ . The inclusion relations between Demazure modules \cite[Proposition 4.4.5]{J} give 
\[
L(\la)_{u \la} \not\subseteq \UU_q^-(\g) L(\la)_{u'' \la}
\]
from which it follows that $(\De_{u \la, \la} +I_w^+) \De_{w \la, \la}^{-1} \in Q^-_w(u'')$. Thus,
\[
d_{u, \la} = \phi_w( (\De_{u \la, \la} +I_w^+) \De_{w \la, \la}^{-1}) \in E_w(u) \cap I_w(u'').
\]
This proves that $E_w(u)$ satisfies the second property for a separating Ore set.
\end{proof}
\subsection{Properties of the constructed separating Ore sets}
\label{3.4}
The $H$-primes and the constructed separating Ore sets of the algebras $R_q[G]$ and $\UU^-[w]$ in Theorems 
\ref{tseparateOre} and \ref{tSetOreUw} have the following properties:

(1) For all $I \in \HSpec R$, $I$ is completely prime.

(2) For $I \in \HSpec R$, there exists an Ore set $E_I$ in $R$, such that $E_I \cap I = \varnothing$, 
the elements $E_I$ are normal modulo $I$ and the localization $R_I := (R/I)[E_I^{-1}]$ is $H$-simple.

(3) For each incident pair $I \subseteq J$, the image of $E_J$ in $R/I$ is a separating Ore set for the pair. 
(The elements of this image are not normal in general.)

\end{document}